\documentclass[a4paper,leqno,12pt]{amsart}
\usepackage{graphicx} 
\usepackage{graphicx,epsfig}
\usepackage{amscd}
\usepackage{amsmath,latexsym,amssymb,amsthm}
\usepackage[all]{xy}
\usepackage[nospace,noadjust]{cite}
\usepackage{setspace,cite}
\usepackage{lscape,fancyhdr,fancybox}
\usepackage{graphicx}

\usepackage{color}
\usepackage{url}
\usepackage{enumerate}
\usepackage[mathscr]{euscript}
\usepackage{tikz-cd}
\input xy
\xyoption{all}

  \theoremstyle{plain}

    \newtheorem{thm}{Theorem}[section]
    \newtheorem{prop}[thm]{Proposition}

    \newtheorem{subsec}[thm]{}

\theoremstyle{definition}
    \newtheorem{defn}[thm]{Definition}
\theoremstyle{remark}

    \newtheorem*{ack}{Acknowledgements}

\newcommand{\Tr}{\operatorname{Tr}}

\newcommand{\R}{\mathbb{R}}
\newcommand{\Q}{\mathbb{Q}}
\newcommand{\C}{\mathbb{C}}

\title{}
\author{}
\date{}
\usepackage{amssymb}

\begin{document}
\title{Endomorphisms of the Cohomology Algebra of the Even Orthogonal Grassmannian}
\author{Arnab Goswami}
\email{arnab.goswami@cbs.ac.in ; arnab543@gmail.com}
\address{School of Mathematical Sciences, 
 UM-DAE Centre for Excellence in Basic Sciences, 
University of Mumbai,  Kalina,  Santacruz (East), 
Mumbai - 400098,  India.}

\author{Swagata Sarkar}
\email{swagata.sarkar@cbs.ac.in ;  swagatasar@gmail.com}
\address{School of Mathematical Sciences, 
 UM-DAE Centre for Excellence in Basic Sciences, 
University of Mumbai,  Kalina,  Santacruz (East), 
Mumbai - 400098,  India.}

\subjclass[2010]{Primary: 57T15, 14M17,  55M99;\ Secondary: 14M15, 55S99}
\keywords{Cohomology Endomorphisms,  Homogeneous Spaces,  Spaces $G/P$}

\maketitle

\begin{abstract}
Let $M_{n,k}$ denote the even orthogonal Grassmanian,  $SO(2n) / (U(k) \times SO(2n-2k) )$.  We study  endomorphisms of the rational cohomology algebra of $M_{n,k}$.  We prove that an endomorphism of the rational cohomology algebra of $M_{n,k}$,  which maps all the Chern classes of the canonical $k$-plane bundle over $M_{n,k}$ to zero,  or maps all the Pontrjagin classes of the canonical,  real,  oriented $(2n-2k)$-plane bundle over $M_{n,k}$ to zero,   is the zero endomorphism.   Additionally,  we prove that if an endomorphism of the rational cohomology algebra of $M_{n,k}$ vanishes on $H^{2}(M_{n,k}; \Q)$,   and admits a splitting,  then the splitting equals zero.  

\end{abstract}

\section{Introduction}

\noindent
We study the endomorphisms of the rational cohomology algebra of $M_{n,k}$ ,  where $M_{n,k}$ is the homogeneous space $SO(2n) / U(k) \times SO(2n-2k) $.  $M_{n,k}$, known as the even orthogonal Grassmannian,  is a space of the type $G/P$,  where $G$ is the complex Lie Group $SO(2n; \C)$ and $P$ is a maximal parabolic subgroup.  It is a compact Hermititan symmetric space and hence,  a Kahler manifold.   \\
\noindent
There is a canonical $k$-plane bundle $\omega$ and a canonical,  real, oriented $(2n-2k)$-plane bundle $\xi$,  such that certain cohomology classes of these bundles generate the rational cohomology algebra of $M_{n,k}$.  We will describe the bundles $\omega$ and $\xi$ in \S 2.  \\
 
 \noindent
We study the properties of certain specific types of endomorphisms of $M_{n,k}$ and prove: 

\begin{thm} \label{prop61}
Let $h \colon H^{*} (M_{n,k} ; \Q) \rightarrow H^{*} (M_{n,k} ; \Q)$ be an endomorphism of the cohomology algebra   $H^{*} (M_{n,k} ; \Q)$,  which maps all Chern classes of $\omega$ to zero (that is,  $h(c_{i}) = 0$ for all $i \in \{1,  \cdots ,  k \}$).  Then $h$ is the zero endomorphism.  
\end{thm}

\noindent
We have an analogous result for the case where $h(p_{j}) = 0$ for all $1 \leq j \leq n-k$. 
 \begin{thm} \label{prop62}
Let $h \colon H^{*} (M_{n,k} ; \Q) \rightarrow H^{*} (M_{n,k} ; \Q)$ be an endomorphism of the cohomology algebra   $H^{*} (M_{n,k} ; \Q)$,  which maps all Pontrjagin classes of $\xi$ to zero (that is,  $h(p_{j}) = 0$ for all $j \in \{1,  \cdots ,  n-k \}$).  Then $h$ is the zero endomorphism.  
\end{thm}

\noindent
\begin{defn}
Let $p^{*} \colon H^{*}(M_{n,k} ; \Q) \rightarrow H^{*}( SO(2n)/ T^{n}; \Q)$ be the monomorphism induced by the map $p \colon SO(2n)/ T^{n}  \longrightarrow  M_{n,k}$,  where $T^{n}$ is a maximal torus of $SO(2n)$.  \\
Then a {\it splitting} of an endomorphism $h: H^{*}( M_{n,k} ; \Q) \to H^{*}(M_{n,k} ; \Q)$,  is an endomorphism $\widetilde{h}: H^{*}( SO(2n)/ T^{n}; \Q) \longrightarrow H^{*}( SO(2n)/ T^{n}; \Q)$,  such that the following diagram commutes: \\

\begin{tikzcd}
H^{*}(M_{n,k},\mathbb{Q})\arrow[r,"h"] \arrow[d,"p^{*}"] & H^{*}(M_{n,k},\mathbb{Q}) \arrow[d ,"p^{*}"] \\ 
 H^{*}(SO(2n)/T^{n}; \mathbb{Q}) \arrow[r,"\widetilde{h}"] & H^{*}(SO(2n)/T^{n}; \mathbb{Q})
\end{tikzcd}
\end{defn}

\noindent
 We prove the following: 
 
 \begin{thm} \label{prop63}
Let $h:H^{*}(M_{n,k}; \Q) \to H^{*}(M_{n,k}; \Q)$  be an endomorphism of the cohomology algebra,  such that $h$ vanishes on $H^{2}(M_{n,k} ,  \Q)$.   If $h$ admits a splitting $\widetilde{h}$,   then $\widetilde{h} =0$,  and $h =0$. 
\end{thm}

\noindent 
We have analogous results when the image of $p_{1}$  is zero ,   and also when $n-k=2$ and the image of $e_{2n-2k}$ is zero.  \\

\begin{thm} \label{prop64}
Let $h:H^{*}(M_{n,k}; \Q) \to H^{*}(M_{n,k}; \Q)$  be an endomorphism of the cohomology algebra,  such that $h(p_{1})=0$,  where $p_{1}$ denotes the first Pontrjagin class of $\xi$.  If $h$ admits a splitting $\widetilde{h}$,   then $\widetilde{h} =0$.  In such a case,  $h$ will be the zero endomorphism of $H^{*}(M_{n,k};\Q)$.\\
\end{thm}

\begin{thm}\label{prop65}
Let $h:H^{*}(M_{n,k}; \Q) \to H^{*}(M_{n,k}; \Q)$  be an endomorphism  such that $h(e_{2n-2k})=0$,  where $n-k=2$ and $e_{2n-2k}$ denotes the Euler class of $\xi$.   If $h$ admits a splitting $\widetilde{h}$,   then $\widetilde{h} =0$.  and $h$ is also the zero endomorphism. 
\end{thm}

\noindent
Note that given any non-zero $\alpha\in \Q$, one has automorphisms $h_{\alpha} \colon H^r(M_{n,k}; \Q) \rightarrow H^r(M_{n,k}; \Q)$ defined as $h_{\alpha}(x)=\alpha^r x$ for all $x\in H^r(M_{n,k}; \Q)$. These are known as the Adams maps,  and satisfy the hypothesis. \\

\noindent
We were motivated by the paper of M. Hoffman \cite{hoff1}, who studied the same problem for complex Grassmann manifolds.  
 See also the works of Brewster (\cite{brew}),  Glover,  Homer,  ({\cite{hogr1},  \cite{hogr2},  \cite{hoho1}),  Papadima (\cite{padi}),  H. Duan (\cite{duan1},  \cite{duan2},  \cite{duan3},  \cite{duan4}),   Shiga-Tezuka  (\cite{shte}),   Lin (\cite{lin}), and Kaji and Theriault (\cite{kath}).  \\

\noindent
The paper is arranged as follows.  In \S 2 we describe the homogeneous space $M_{n,k}$ and compute its rational cohomology.   In \S 3,  we give the proofs of our main results. In \S 4 we discuss the Lefschetz number of certain maps.  \\

\begin{ack}
The authors would like to acknowledge and thank Professors  Samik Basu,  R.  C.  Cowsik,  M.  S.  Raghunathan,  and Parameswaran Sankaran for useful mathematical discussions,  and their many helpful comments and suggestions during  the writing of this paper.  \\
The authors are grateful to H. Duan for pointing out a gap in one of the proofs in a previous version of our paper.  
\end{ack}

\vspace{0.3cm}
 
\section{The Homogeneous Space $M_{n,k}$ and its Cohomology}

\noindent
 We consider the even orthogonal Grassmannian,  $M_{n,k} = SO(2n) / (U(k) \times SO(2n-2k))$.  
This quotient space, $M_{n,k}$, is a homogeneous space of the form $G/P_{k}$,  where $G$ is a simple,  complex,  connected Lie group of the type $D_{n}$ and $P_{k}$,  a maximal,  parabolic subgroup of $G$.  $M_{n,k}$ is a compact,  complex manifold of real dimension $2(2nk - k^{2}) -k(k+1)$.  When $k=n$,  this space is the Grassmannian of complex structures,  $SO(2n)/U(n)$.  In this paper,  we study the  endomorphisms of the rational cohomology algebra of $M_{n,k}$.  \\

\noindent
Next we compute the cohomology of the space $M_{n,k}$.  
For a compact,  connected Lie Group $G$,  let $T$ denote a maximal torus,  $W(G)$  denote its Weyl group,  and $BG$ denote the classifying space of $G$.  \\

\noindent
Consider the universal complex $k$-plane bundle $\gamma_{k}^{\C}$ over $BU(k)$ and the universal real $(2n-2k)$-bundle $\gamma_{2n-2k}^{\R}$ over $BSO(2n-2k)$.  Let $ c_{i} (\gamma_{k}^{\C}) \in H^{2i} (BU(k) ; \Q)$ be the $i$-th Chern class of $\gamma_{k}^{\C}$ over $BU(k)$.  Let $ p_{j} (\gamma_{2n-2k}^{\R}) \in H^{2j} (BSO(2n-2k) ; \Q)$ be the $j$-th Pontrjagin class of $\gamma_{2n-2k}^{\R}$ over $BSO(2n-2k)$. Let $e_{2n-2k} (\gamma_{2n-2k}^{\R}) \in H^{2n-2k}(BSO(2n-2k) ; \Q)$ be the Euler class of $\gamma_{2n-2k}^{\R}$ over $BSO(2n-2k)$.  Note that $p_{n-k} (\gamma_{2n-2k}^{\R}) = e_{2n-2k}^{2} (\gamma_{2n-2k}^{\R})$. \\

\noindent
Let $G = SO(2n)$  and $H = U(k) \times SO(2n-2k)$.  Consider  $\omega$,  the canonical complex $k$-plane vector bundle $G \times_{H} \C^{k} \longrightarrow G/H = M_{n,k}$ and $\xi$,  the canonical real,  oriented $(2n-2k)$-plane vector bundle $G \times_{H} \R^{2n-2k} \longrightarrow G/H = M_{n,k}$.  Then there are classifying maps $f_{\omega} \colon M_{n,k} \longrightarrow BU(k)$ and $f_{\xi} \colon M_{n,k} \longrightarrow BSO(2n-2k)$.  Therefore,  we have a map $f \colon M_{n,k} \longrightarrow BU(k) \times BSO(2n-2k)$ such that the pullbacks $f^{*} (\gamma_{k}^{\C}) = \omega$ and $f^{*} (\gamma_{2n-2k}^{\R}) = \xi$.  Note that $f$ induces a surjection in rational cohomology.  \\

\noindent
Then we have the following description of the rational cohomology algebra of $M_{n,k}$:
\begin{thm} \label{thm31}
The rational cohomology algebra  $H^{*}(M_{n,k} ; \Q)$ is of the form $$\Q[c_1 , \cdots , c_k] \otimes \Q[p_1,  \cdots , p_{n-k-1},  e_{2n-2k}]/ I_{n,k} $$ where $c_{i}$'s are the Chern classes of the bundle $\omega$,  $p_j$'s are the Pontrjagin classes of the bundle $\xi$ and $e_{2n-2k}$ is the Euler class of the bundle $\xi$,  and $I_{n,k}$ is the ideal generated by the polynomials with constant term equal to zero,  which are invariant under the action of the Weyl group of $SO(2n)$. 
\end{thm}

\begin{proof}
\noindent
Recall that (\cite{MMT}),   $H^{*}(BSO(2n); \Q) =  \Q[p_{1},  \cdots ,  p_{n-1},  e_{2n}]$,  where $p_{i}$'s denote the  Pontrjagin classes,  and $e_{2n}$ denotes the Euler class of the universal oriented $2n$-plane bundle over $BSO(2n)$.  (Note that $e_{2n}^{2} = p_{n}$).  The above is a polynomial ring,  and $H^{odd} = 0 $ ,  while $H^{even}$ is a free  module.  
Also recall that  $H^{*}(BU(k); \Q) =  \Q[c_{1},  \cdots ,  c_{k}]$,  where $c_{j}$'s denote the Chern classes of the universal $k$-plane bundle over $BU(k)$.  So we have (\cite{MMT}) 
\begin{displaymath}
H^{*}(BU(k) \times BSO(2n-2k) ; \Q) = H^{*}(BU(k); \Q) \otimes H^{*}( BSO(2n-2k) ; \Q) 
\end{displaymath}
 \begin{displaymath}
= \Q[c_{1}, \cdots, c_{k}] \otimes \Q[p_{1},  \cdots,  p_{n-k-1},  e_{2n-2k}]
\end{displaymath}

\noindent
Note that (\cite{MMT}) for a maximal torus $T^{n}$,  $H^{*}(SO(2n)/T^{n} ; \Q)$  is a finitely generated,  free module and $H^{odd}(SO(2n)/T^{n} ; \Q)=0 $.  In fact,  $H^{*}(SO(2n)/T^{n} ; \Q) \cong \Q[t_{1}, \cdots ,  t_{n}]/I $,  where for each $i$,  $\deg(t_i) =2$ and $I$ is the ideal generated by positive degree polynomials which are invariant under the action of the Weyl group of $SO(2n)$.  \\
 
 \noindent
Now consider the fibration,  \\
\begin{tikzcd}
(U(k)\times SO(2n-2k))/ T^{n} \arrow[r,"j"]  & SO(2n)/ T^{n}  \arrow[r,"p"] & M_{n,k}
\end{tikzcd}  

\noindent
By the Leray-Hirsch Theorems (\cite{MMT}),  we have that $j^{*}$ is an epimorphism and $p^{*}$ is a monomorphism.  Since $H^{*}(M_{n,k} ;\Q)$  maps injectively into $H^{*}(SO(2n)/ T^{n} ; \Q)$ , $H^{*}( M_{n,k} ; \Q)$ is also free and  $H^{odd} (M_{n,k} ; \Q) = 0$. \\

\noindent
Since $H^{odd} (BSO(2n)) = 0$ and $H^{odd} (M_{n,k} ; \Q) = 0$,  the differential $d_{r}^{p,q}$ on each page of the spectral sequence is zero,  and, we have that the Serre Spectral sequence of the fibration 

\begin{tikzcd}
 M_{n,k} \arrow[r,"j_{1}"]  & B(U(k) \times SO(2n-2k)) \arrow[r,"q"] & BSO(2n)
\end{tikzcd} 
\\
 collapses.  Again,  by the Leray-Hirsch Theorem,  we have a surjective map,  as follows: 

\begin{tikzcd}
 H^{*} (B(U(k) \times SO(2n-2k)); \Q )  \arrow[r,"j_{1}^{*}"] & H^{*} (M_{n,k} ; \Q ) 
\end{tikzcd}

\noindent
Let $I_{n,k}$ be the ideal generated by the positive degree terms in the image of $q^{*}$.  From the work of Borel \cite{bor},  it is known that there is an isomorphism:

 \begin{tikzcd}
H^{*}(M_{n,k} ; \Q) \arrow[r,"\simeq"]  & H^{*}(BU(k) \times BSO(2n-2k) ; \Q) /  I_{n,k}
 \end{tikzcd}

\noindent
So, the cohomology ring of $M_{n,k}$ is of the form 
\begin{displaymath}
 H^{*}(M_{n,k}; \Q ) \simeq( \Q [c_{1},...,c_{k}] \otimes \Q [p_{1}, \cdots ,p_{n-k-1},e_{2n-2k}]) /I_{n,k}
\end{displaymath}
\noindent
 where,  $c_{i}$'s are the Chern classes of $\omega$,  $p_{j}$'s are the Pontrjagin classes of $\xi$ and $e_{2n-2k}$ is the Euler class of $\xi$.
  \noindent
 Here,  $I_{n,k} := \langle Im (q^{*}) \rangle $,  is the ideal generated by polynomials with constant term equal to zero,  which are invariant under the action of the Weyl group of $SO(2n)$.  \\
 \noindent
The generating relations of $I_{n,k}$ are obtained as follows.  We have,  $\omega_{\R} \oplus \xi \cong 2n \epsilon_{\R}$,  where $\epsilon$ is the trivial bundle.  Complexifying,  we get,  $\omega \otimes \C \oplus \xi \otimes \C \cong 2n \epsilon_{\C}$.  Hence,  $c(\omega \otimes \C) \cdot c(\xi \otimes \C) =1$,  where $c$ denotes the total Chern class.  \\ Additionally,  $\omega_{\R} \oplus \xi \cong 2n \epsilon_{\R} $  implies  $p(\omega_\R)\cdot  p(\xi) = 1$,  where $p$ denotes the total Pontrjagin class. 
It also implies $ e(\omega_{\R} \oplus \xi) = 0$,  which implies $e(\omega_{\R}) \cdot e(\xi) = 0$,  which,  in turn, implies $c_k \cdot e_{2n-2k} =0$,  since $c_k = e(\omega_{\R})$.  
 
 \noindent
 Note that,  since $p^{*} \colon H^{*} (M_{n,k}) \rightarrow H^{*}(SO(2n)/T^{n})$ is a monomorphism,  and $H^{*}(M_{n,k} ; \Q)$ sits injectively in $\Q[t_{1},  \cdots , t_{n}]/  I  $.  Specifically,  the Chern classes and the Pontrjagin classes can be expressed in terms of symmetric polynomials,  with $c_i = \sigma_i (t_1,  \cdots ,  t_k) $,  
$ p_j = \sigma_j (t_{k+1}^{2},  \cdots , t_{n}^{2})$,  for $1 \leq i \leq k$ and $1 \leq j \leq n-k$,  and $e_{2n-2k} = t_{k+1} \cdots t_{n}$.

\end{proof}

\vspace{0.5cm}

\section{Main Results}

\noindent
This section contains the proofs of our main results,  about the properties of specific endomorphisms of the rational cohomology algebra of the space $M_{n,k}$.  \\

\newtheorem*{myprop}{Theorem}
\begin{myprop} {\bf [Theorem \ref{prop61}]} 
Let $h \colon H^{*} (M_{n,k} ; \Q) \rightarrow H^{*} (M_{n,k} ; \Q)$ be an endomorphism of the cohomology algebra   $H^{*} (M_{n,k} ; \Q)$,  which maps all Chern classes of $\omega$ to zero (that is,  $h(c_{i}) = 0$ for all $i \in \{1,  \cdots ,  k \}$).  Then $h$ is the zero endomorphism.  
\end{myprop}

\begin{proof}
We have the relation from the ideal $I_{n.k}$,  given by $c_{1}^{2} -2c_{2}+p_{1} =0$.  Applying the endomorphism,  $h$,  to this relation,  we get $h(p_{1}) = 0$.  Similarly,  applying $h$ to the general form of the relation in the ideal $I_{n,k}$ : 
\begin{equation}
      (c_{j}^{2}-2c_{j-1}c_{j+1}+...)+(c_{j-1}^{2}-c_{j-2}c_{j}+...)p_{1}+...+(c_{1}^{2}-2c_{2})p_{j-1}+p_{j}=0 \label{eqn67}
  \end{equation}
we get that $h(p_{j}) = 0$  for all $j \in \{1,  \cdots ,  n-k\}$ .  \\
The only relation in $I_{n,k}$,  which contains $e_{2n-2k}$ is the relation:  $c_{k} e_{2n-2k} = 0$.  We now calculate $h(e_{2n-2k} )$.  
Suppose $h(e_{2n-2k}) = \Sigma_{\alpha} A_{\alpha} c_{1}^{\alpha_{1}} \cdots c_{k}^{\alpha_{k}} + Be_{2n-2k}$,  where $\alpha = (\alpha_{1} ,  \cdots , \alpha_{k})$,  $\alpha_{i} \in \Q$,  and $\Sigma_{i=1}^{k} i \alpha_{i} = n-k$.   
Then  $h(p_{n-k}) = h(e_{2n-2k}^{2}) = (\Sigma_{\alpha} A_{\alpha} c_{1}^{\alpha_{1}} \cdots c_{k}^{\alpha_{k}} + B e_{2n-2k})^{2} = 0$.  
We know that $p_{n-k}$ is of degree $4(n-k)$.  \\
\noindent
If for all $i \in \{1,  \cdots k \}$,  $A_{\alpha_{i}} = 0$,  then $B^{2} e_{2n-2k}^{2} =0$,  and hence $B = 0$.  Therefore,  $h(e_{2n-2k}) = 0$.  \\
\noindent
Similarly,   if $B =0$,  then in degree $4(n-k)$,   
\begin{equation}
(\Sigma_{\alpha} A_{\alpha} c_{1}^{\alpha_{1}} \cdots c_{k}^{\alpha_{k}} )^{2} = 0
\label{eqn68}
\end{equation}
\noindent
We have the relation 
\begin{equation}
\label{eqn681}
 (c_{j}^{2} - 2c_{j-1}c_{j+1} + 2c_{j-2}c_{j+2}+ \cdots ) + \cdots + (c_{1}^{2} -2c_{2}) p_{n-k-1} + p_{n-k} =0
 \end{equation}
Now,  $p_{n-k}$ is a function of $c_{1},  \cdots c_{k}$.  Therefore,  
\begin{equation}
f(c_{1},  \cdots ,  c_{k}) + p_{n-k} = 0
\label{eqn69}
\end{equation}
\noindent
Note that both \ref{eqn68} and \ref{eqn69} are relations in degree $4(n-k)$.  
Let,  if possible,  $f(c_{1},  \cdots ,  c_{k}) = (\Sigma_{\alpha} A_{\alpha} c_{1}^{\alpha_{1}} \cdots c_{k}^{\alpha_{k}} )^{2}$.  Then from \ref{eqn68} and \ref{eqn69},  we get $p_{n-k} =0$. This is impossible since Pontrjagin classes are not zero.  Hence $f(c_{1},  \cdots ,  c_{k}) \neq (\Sigma_{\alpha} A_{\alpha} c_{1}^{\alpha_{1}} \cdots c_{k}^{\alpha_{k}} )^{2}$.   But $c_{1}^{\alpha_{1}} \cdots c_{k}^{\alpha_{k}}$ and $p_{n-k}$ are generators in degree $4(n-k)$.  Comparing \ref{eqn68} and \ref{eqn69},  we get that $A_{\alpha}$ equals zero for each $\alpha $.  Therefore,  $h(e_{2n-2k}) = 0$.  \\
\noindent
If not all $A_{\alpha}$'s are equal to zero and $B$ is also not equal to zero,  we have \\
 $h(e_{2n-2k}) = \Sigma_{\alpha} A_{\alpha} c_{1}^{\alpha_{1}} \cdots c_{k}^{\alpha_{k}} + Be_{2n-2k}$.  
 Therefore,  
  \begin{equation}
  \label{eqn682}
  0 = h(p_{n-k}) = h(e_{2n-2k}^{2}) = 
(\Sigma_{\alpha} A_{\alpha} c_{1}^{\alpha_{1}} \cdots c_{k}^{\alpha_{k}} + Be_{2n-2k})^{2}.
\end{equation}
Hence,   
$ (\Sigma_{\alpha}A_{\alpha} c_{1}^{\alpha_{1}} \cdots c_{k}^{\alpha_{k}})^{2} + B^{2} e_{2n-2k}^{2} + 2B (\Sigma_{\alpha} A_{\alpha} c_{1}^{\alpha_{1}} \cdots c_{k}^{\alpha_{k}}) e_{2n-2k} = 0$.  \\
Now degree of $c_{1}^{\alpha_{1}} \cdots c_{k}^{\alpha_{k}} e_{2n-2k}$ is $4(n-k)$.  But $e_{2n-2k}$  only appears in one relation (in degree $2n$): $c_{k} e_{2n-2k} = 0$.  It does not appear in any relation in degree equal to $4(n-k)$.  Therefore,  comparing Equations $\ref{eqn681}$ and $\ref{eqn682}$,  we get,  $A_{\alpha}$ is zero for every $\alpha$.  
Hence,  $h(e_{2n-2k}) = 0$,  and from the above discussion,  $h$ is the zero endomorphism.  

\end{proof}

 \begin{myprop}  {\bf [Theorem \ref{prop62}]} 
 Let $h \colon H^{*} (M_{n,k} ; \Q) \rightarrow H^{*} (M_{n,k} ; \Q)$ be an endomorphism of the cohomology algebra   $H^{*} (M_{n,k} ; \Q)$,  which maps all Pontrjagin classes of $\xi$ to zero (that is,  $h(p_{j}) = 0$ for all $j \in \{1,  \cdots ,  n-k \}$).  Then $h$ is the zero endomorphism.  
\end{myprop}
 
 \begin{proof}
 Since $h(p_{j})= 0$ for each $j \in \{1,  \cdots , n-k \}$,  we have $h(p_{n-k}) = h(e_{2n-2k}^{2}) =0$.  Similar to proof of Theorem \ref{prop61},  we can show that $h(e_{2n-2k}) = 0$.  Now,  applying $h$ to the relations in the ideal,  and using $h(p_{j}) =0$,  we have the following: 
 \begin{eqnarray}
    h(c_{1})^{2}-2h(c_{2}) &=&0  \nonumber \\
    h(c_{2})^{2}- 2h(c_{1})h(c_{3})+2h(c_{4})&=& 0 \nonumber  \\
    h(c_{3})^{2}-2h(c_{2})h(c_{4})+2h(c_{1})h(c_{5})-2h(c_{6}) &=& 0 \nonumber \\
    \cdots  \nonumber \\
    h(c_{k-2})^{2}- 2h(c_{k-3})h(c_{k-1})+ 2h(c_{k-4})h(c_{k}) &=& 0 \nonumber \\
    h(c_{k-1})^{2}-2h(c_{k-2})h(c_{k}) &=& 0 \nonumber \\
    h(c_{k})^{2}&=&0 \nonumber
\end{eqnarray}

\noindent
We show that  $h(c_{k}) =0$.  The proof is divided into three cases.  \\
\noindent
When $k <n-k$,   we can write $k+s = n-k$,  where $s$ is a positive integer. Then $h(c_{k})$ is of the form: 
$h(c_{k}) = \Sigma_{\alpha} A_{\alpha} c_{1}^{\alpha_{1}} \cdots c_{k}^{\alpha_{k}} $,  
where $\alpha = (\alpha_{1} ,  \cdots ,  \alpha_{k})$ ,  $\alpha_{i}$'s are non-negative integers,  $A_{\alpha}$'s are rational numbers,  and $\Sigma_{i=1}^{k} i \alpha_{i} = k$.    We have the relation 
\begin{equation}
 \label{eqn36}
    h(c_{k})^{2}= (\Sigma_{\alpha} A_{\alpha} c_{1}^{\alpha_{1}}\dots c_{k}^{\alpha_{k}})^{2} = 0
\end{equation}
\noindent
Note that this relation is in degree $4k$.  In the  ideal  there exists a relation in degree $4k$,  of the following form: 
\begin{equation}
 \label{eqn37}
    c_{k}^{2}+(c_{k-1}^{2}-2c_{k-2}c_{k})p_{1}+ \cdots +(c_{1}^{2}-2c_{2})p_{n-k-s-1}+p_{n-k-s}=0
\end{equation}

\noindent
Substituting the values of $p_{1},  \cdots,  p_{n-k-s-1}$ in Equation $\ref{eqn37}$,  one sees that $p_{n-k-s}$ is always present in Equation $\ref{eqn37}$ but not in Equation $\ref{eqn36}$.   Comparing equations $\ref{eqn36}$ and $\ref{eqn37}$,  we get $A_{\alpha}=0$ for each $\alpha$.   Therefore,   $h(c_{k})=0$.\\  

\noindent
In the case where $k = n-k$,    $h(c_{k})$ of the form:  
$h(c_{k}) = \Sigma_{\alpha} A_{\alpha} c_{1}^{\alpha_{1}} \cdots c_{k}^{\alpha_{k}} + Be_{2n-2k}$,  
where $\alpha = (\alpha_{1} ,  \cdots ,  \alpha_{k})$ ,  $\alpha_{i}$'s are non-negative integers,  $A_{\alpha}$,  $B$ are rational numbers,  and $\Sigma_{i=1}^{k} i \alpha_{i} = k = n-k$.    Therefore,  $h(c_{k})^{2}= (\Sigma_{\alpha} A_{\alpha} c_{1}^{\alpha_{1}}\dots c_{k}^{\alpha_{k}}+Be_{2n-2k})^{2}=0$. \\
In the  ideal ,  there is a relation in degree $4k$ 
\begin{equation}
\label{eqn38}
    c_{k}^{2}+(c_{k-1}^{2}-2c_{k-2}c_{k})p_{1}+...+(c_{1}^{2}-2c_{2})p_{n-k-1}+p_{n-k}=0
\end{equation}
\noindent
We have 
 \begin{equation}
 \label{eqn39}
    h(c_{k})^{2}= (\Sigma_{\alpha} A_{\alpha} c_{1}^{\alpha_{1}}\dots c_{k}^{\alpha_{k}})^{2}+B^{2}p_{n-k}+2B(\Sigma_{\alpha} A_{\alpha} c_{1}^{\alpha_{1}}\dots c_{k}^{\alpha_{k}})e_{2n-2k} = 0
 \end{equation}
\noindent
Comparing  Equations $\ref{eqn38}$ and $\ref{eqn39}$ , we get $A_{\alpha} \cdot B=0$,  except when $\alpha=(0,...,0,1)$,  since $c_{k}e_{2n-2k}=0$. \\

\noindent
If $B=0$ then $h(c_{k})^{2}=(\Sigma_{\alpha} A_{\alpha} c_{1}^{\alpha_{1}}\dots c_{k}^{\alpha_{k}})^{2} = 0$.  By a similar argument to the previous case,  we have  $h(c_{k})=0$.\\
  \noindent
If $B\neq 0$,  $h$ is of this following form: $h(c_{k})=A_{k}c_{k}+Be_{2n-2k}$
Comparing $h(c_{k})^{2}=A_{k}^{2}c_{k}^{2}+B^{2}p_{n-k} = 0$ with the relation \ref{eqn38},  we have $A_{k}=0=B$.  Therefore,  $h(c_{k})=0$.\\

\noindent
Now consider the last case,  $k>n-k$.  In this case,   $k-s = n-k$,  where $0<s<k$,  and $h(c_{k})$ is of the form: 
 $h(c_{k}) = \Sigma_{\alpha} A_{\alpha} c_{1}^{\alpha_{1}} \cdots c_{k}^{\alpha_{k}} + (\Sigma B_{\beta} c_{1}^{\beta_{1}} \cdots c_{s}^{\beta_{s}})e_{2n-2k}$, 
where $\alpha = (\alpha_{1} ,  \cdots ,  \alpha_{k})$ ,  $\beta = (\beta_{1} ,  \cdots ,  \beta_{s})$, $\alpha_{i}$'s and $\beta_{j}$'s are non-negative integers,   $\Sigma_{i=1}^{k} i \alpha_{i} = k = n-k+s$,  and $\Sigma_{i=1}^{s} i \beta_{i} = s$.   \\

\noindent
If all $A_{\alpha} = 0$,  then  $h(c_{k}) =  (\Sigma B_{\beta} c_{1}^{\beta_{1}} \cdots c_{s}^{\beta_{s}})e_{2n-2k}$.
Therefore,  $h(c_{k}^{2}) = (\Sigma B_{\beta} c_{1}^{\beta_{1}} \cdots c_{s}^{\beta_{s}})^{2} p_{n-k}$.  If we have $h(c_{k}^{2}) =0$,  then $ (\Sigma B_{\beta} c_{1}^{\beta_{1}} \cdots c_{s}^{\beta_{s}})^{2} p_{n-k} = 0$. Now,  $p_{n-k}$ is a function of $c_{1} ,  \cdots c_{k}$.  We can write  $f(c_{1} ,  \cdots ,  c_{k}) + p_{n-k} =0$.  
Therefore,  
\begin{equation}
 (\Sigma B_{\beta} c_{1}^{\beta_{1}} \cdots c_{s}^{\beta_{s}})^{2} f(c_{1} ,  \cdots ,  c_{k}) = 0.   
 \label{eqn610}
 \end{equation}
 This is a relation of degree $4k$.  In degree $4k$,  the relation is given by 
\begin{equation}
c_{k}^{2} + (c_{k-1}^{2} -c_{k-2}c_{k})p_{1} + \cdots + (c_{s}^{2} - 2c_{s-1}c_{s+1} + \cdots ) p_{n-k} = 0.  
 \label{eqn611}
\end{equation}
Comparing the above two relations (\ref{eqn610}) and (\ref{eqn611}) we observe that the term $c_{k}^{2}$,  which appears in (\ref{eqn611}),  never appears in (\ref{eqn610}).  Therefore,  $B_{\beta} = 0$ for every $\beta$,  and $h(c_{k}) =0$. \\
\noindent
If, at least one of $A_{\alpha}$'s and $B_{\beta}$'s is non-zero,  then 
$h(c_{k}) = \Sigma_{\alpha} A_{\alpha} c_{1}^{\alpha_{1}} \cdots c_{k}^{\alpha_{k}} + (\Sigma B_{\beta} c_{1}^{\beta_{1}} \cdots c_{s}^{\beta_{s}})e_{2n-2k}$
\noindent
Since $h(c_{k}^{2}) = 0$,  we get 
\begin{equation}
 h(c_{k})^{2}= ( \Sigma_{\alpha} A_{\alpha} c_{1}^{\alpha_{1}} \cdots c_{k}^{\alpha_{k}} + (\Sigma B_{\beta} c_{1}^{\beta_{1}} \cdots c_{s}^{\beta_{s}})e_{2n-2k})^{2} = 0  
\label{eqn612}
\end{equation}
\begin{equation}
=  (\Sigma_{\alpha} A_{\alpha} c_{1}^{\alpha_{1}}c_{2}^{\alpha_{2}}...c_{k}^{\alpha_{k}})^{2}+ (\Sigma_{\beta} B_{\beta} c_{1}^{\beta_{1}}...c_{s}^{\beta_{s}})^{2}p_{n-k} + 2\Sigma_{\alpha, \beta}A_{\alpha}B_{\beta}(c_{1}^{\alpha_{1}}c_{2}^{\alpha_{2}}...c_{k}^{\alpha_{k}})(c_{1}^{\beta_{1}}...c_{s}^{\beta_{s}})e_{2n-2k}
\nonumber
\end{equation}
\noindent
Now,  comparing relations (\ref{eqn611}) and (\ref{eqn612}),  (both of which are in degree $4k$),  we find that (\ref{eqn611}) does not have a term of the form $(c_{1}^{\alpha_{1}} \cdots c_{k}^{\alpha_{k}}) (c_{1}^{\beta_{1}} \cdots c_{s}^{\beta_{s}})e_{2n-2k}$.  Hence,  $A_{\alpha} =0$ for each $\alpha$ and $B_{\beta}=0$ for each $\beta$.  Therefore $h(c_{k}) =0$.  \\
\noindent
Next we consider the case where $B_{\beta} =0$,  for all $\beta$.  Then $h(c_{k}) = \Sigma_{\alpha} A_{\alpha} c_{1}^{\alpha_{1}} \cdots c_{k}^{\alpha_{k}} $.  There are two subcases to be considered: \\
\noindent
{\it Subcase 1:} $ k-(n-k) \geq 2$.  Assume $ h(c_{k})^{2}= 0$ 
We have , 
\begin{equation}
    h(c_{k}^{2}) = (\Sigma_{\alpha} A_{\alpha} c_{1}^{\alpha_{1}}c_{2}^{\alpha_{2}}...c_{k}^{\alpha_{k}})^{2} =0
\label{eqn614}
\end{equation}
\noindent
This is a relation in degree $4k$.  We also have the following relation in degree $4k$:
\begin{equation}
    c_{k}^{2}+(c_{k-1}^{2}-2c_{k-2}c_{k})p_{1}+(c_{k-2}^{2}-2c_{k-3}c_{k-1}+2c_{k-4}c_{k})p_{2}+...+(c_{2}^{2}-2c_{1}c_{3}+2c_{4})p_{n-k} =0
    \nonumber
\end{equation}
\noindent
Substituting the values of $p_{1}, \cdots ,p_{n-k}$ in (\ref{eqn611}),  we  observe that there is a term $c_{k-1}^{2}c_{2}$.
\begin{equation}
  c_{k}^{2}+(c_{k-1}^{2}-2c_{k-2}c_{k})(-c_{1}^{2}+2c_{2})+ \cdots = 0  
\end{equation}
\noindent
The coefficient of $c_{k-1}^{2} c_{2}$ is $(2+4m)$,  for some integer $m$ .  Now $2+4m \neq 0$,  since $m$ is an integer.  Comparing the relations in degree $4k$,  we see that $c_{k-1}^{2}c_{2}$ never appears in \ref{eqn614}.  Hence $A_{\alpha}=0$  for each $\alpha$,  and  $h(c_{k})=0$.\\
\noindent
{\it Subcase 2:} $n-k=k-1$.  
We have the following relation in degree $4k$: 
\begin{equation}
  c_{k}^{2}+(c_{k-1}^{2}-2c_{k-2}c_{k})p_{1}+ \cdots +(c_{2}^{2}-2c_{1}c_{3}+2c_{4})p_{k-2}+ (c_{1}^{2}-2c_{2})p_{k-1} =0
  \label{eqn615}
\end{equation}
\noindent
We know that $h(c_{k})$ is a linear combination (over rationals) of $c_{k}, c_{k-1}c_{1}, c_{1}^{k} ,  \cdots$ (the linear combinations vary over all partitions of $k$): 
\begin{equation}
    h(c_{k}) = A_{{1}^{k}}c_{1}^{k}+ \cdots +A_{{k-11}}c_{k-1}c_{1}+ \cdots +A_{k}c_{k}
    \nonumber
\end{equation}
Since $h(c_{k})^{2}=0$, we get 
\begin{equation}
    0 = h(c_{k})^{2}= (A_{{1}^{k}}c_{1}^{k}+ \cdots +A_{{k-11}}c_{k-1}c_{1}+ \cdots +A_{k}c_{k})^{2}
     \label{eqn616}
\end{equation}
\noindent
Substituting values of $p_{1}, \cdots ,p_{k-1}$ in (\ref{eqn615}) ,  one observes that the coefficient of $c_{k-1}^{2}c_{1}^{2}$ is $(-2+4 m)$,  for some integer $m$,  which is not equal to zero. \\
(\ref{eqn615}) is a relation in the ideal,  of degree $4k$ ,  (\ref{eqn616}) is also an relation of the degree $4k$.\\
The term $c_{k-1}^{2} c_{1}^{2}$ always appears in (\ref{eqn615}),  but if $A_{k-1,1} = 0$,  $c_{k-1}^{2} c_{1}^{2}$ never appears in $h(c_{k}^2)$.  Therefore,  comparing $h(c_{k}^2) =0$ and (\ref{eqn615}),  we get that each of  $A_{1^{k}},  \cdots ,  A_{k}$ equals zero.  
Hence $h(c_{k}) = 0$.  \\
In case $A_{k-1,1} \neq 0$,  in $h(c_{k})$ then 
\begin{equation}
c_{k}^{2} + (-2+4m) c_{k-1}^{2}c_{1}^{2} + \cdots = 0 
\label{eqn618}
\end{equation}
 Taking $h(c_{k}^{2}) = 0$,  and substituting the value of $c_{k-1}^{2}c_{k}^{2}$ from \ref{eqn618},  we get 
 \begin{equation}
 A_{k}^{2}c_{k}^{2} + \cdots + \frac{A_{k-1,1}}{2-4m} (c_{k}^{2} + \cdots ) + \cdots =0
 \end{equation}
Therefore,  $\dfrac{A_{k}^{2}}{A_{k-1,1}^{2}} = \dfrac{-1}{2-4m}$.  \\
Since $A_{k}$ and $A_{k-1,1}$ are rational,  $\dfrac{2A_{k}}{A_{k-1,1}}$ is rational,  but $\left ( \dfrac{2}{2m-1} \right )^{\frac{1}{2}}$ is not rational.  Therefore,  $A_{k-1,1} =0$,  and hence,  $A_{\alpha} =0$ for all $\alpha$,  and $h(c_{k}) =0$.  Since we know that $h(c_{k}^{2})=0$,  we have proved that $h(c_{k}) = 0$.  \\
\noindent
Now consider the relation $h(c_{k-1})^{2} - h(c_{k-2}) h(c_{k}) =0$.  Since $h(c_{k}) = 0$,  $h(c_{k-1})^{2} =0$.   In general,  $h(c_{k-1})$ is of the form 
 $h(c_{k-1}) = \Sigma_{\alpha} A_{\alpha} c_{1}^{\alpha_{1}} \cdots c_{k}^{\alpha_{k}} + (\Sigma B_{\beta} c_{1}^{\beta_{1}} \cdots c_{s}^{\beta_{s}})e_{2n-2k}$,  \\
where $\alpha = (\alpha_{1} ,  \cdots ,  \alpha_{k})$ ,  $\beta = (\beta_{1} ,  \cdots ,  \beta_{s})$, $\alpha_{i}$'s and $\beta_{j}$'s are non-negative integers,   $\Sigma_{i=1}^{k} i \alpha_{i} = k-1$,  and $\Sigma_{j=1}^{s} j \beta_{j} =(k-1) - (n-k)$.  \\
Now $h(c_{k-1}^{2})=0$ implies that 
\begin{equation} 
(\Sigma_{\alpha} A_{\alpha} c_{1}^{\alpha_{1}} \cdots c_{k}^{\alpha_{k}} + (\Sigma B_{\beta} c_{1}^{\beta_{1}} \cdots c_{s}^{\beta_{s}})e_{2n-2k})^{2} =0.
\label{eqn620}
\end{equation} 
This is a relation in degree $4(k-1)$.  Consider the relation (also in degree $4(k-1)$): 
\begin{equation}
    (c_{k-1}^{2}-2c_{k-2}c_{k})+(c_{k-2}^{2}-2c_{k-3}c_{k-1}+2c_{k-4}c_{k})p_{1}+ \cdots + (c_{1}^{2}-2c_{2})p_{n-k}=0
    \label{eqn619}
\end{equation}
The term $c_{k-2} c_{k}$ appears in \ref{eqn619},  but does not appear in \ref{eqn620}.  Therefore,  $A_{\alpha} =0$ for each $\alpha$ and $B_{\beta} =0$ for each $\beta$.  Hence,  $h(c_{k-1}) =0$.  \\

\noindent
Now,  using 
\begin{equation}
h(c_{i})^{2} = (\Sigma_{\alpha} A_{\alpha} c_{1}^{\alpha_{1}} \cdots c_{i}^{\alpha_{i}} + (\Sigma B_{\beta} c_{1}^{\beta_{1}} \cdots c_{s}^{\beta_{s}})e_{2n-2k})^{2}
\end{equation}
we get that $A_{\alpha} =0$ for each $\alpha$ and $B_{\beta} =0$ for each $\beta$.  This is because  the term $2c_{i-1} c_{i+1}$ never appears in $h(c_{i})^{2}$.  
Hence,  for each $i$,  $1\leq i \leq k$,  we have $h(c_{i}) =0$.   \\
\noindent
Therefore,  $h$ is the zero endomorphism,  and we are done. 

 \end{proof}

 \begin{myprop}{\bf [Theorem \ref{prop63}]} 
 Let $h:H^{*}(M_{n,k}; \Q) \to H^{*}(M_{n,k}; \Q)$  be an endomorphism of the cohomology algebra,  such that $h$ vanishes on $H^{2}(M_{n,k} ,  \Q)$.   If $h$ admits a splitting $\widetilde{h}$,   then $\widetilde{h} =0$,  and $h =0$. 
\end{myprop}

\begin{proof}
We have the following commuting diagram: 

 \begin{tikzcd}
H^{*}(M_{n,k},\mathbb{Q})\arrow[r,"h"] \arrow[d,"p^{*}"] & H^{*}(M_{n,k},\mathbb{Q}) \arrow[d ,"p^{*}"] \\ 
 H^{*}(SO(2n)/T^{n}; \mathbb{Q}) \arrow[r,"\widetilde{h}"] & H^{*}(SO(2n)/T^{n}; \mathbb{Q})
\end{tikzcd}

\noindent
Now $h(c_{1})=0$ implies $p^{*}\circ h(c_{1})= \widetilde{h}\circ p^{*}(c_{1}) = 0$.  That is,   $ \widetilde{h}(t_{1}+\cdots +t_{k})=0$.  \\
Let $ H = (a_{ij})_{1 \leq i, j \leq n}$ be the matrix representation of the endomorphism $\widetilde{h}$.  
From  $H$,  we get,   $a_{1i}+a_{2i}+...+a_{ki}=0$,  for all $i \in \{1,2,...,n\}$.  
This means that the first $k$ rows of the matrix $H$ are linearly dependent,  and hence,  $H$ is a singular matrix.  
For $i \in \{1,2,...,n\}$,  let the columns of the matrix $H$ be denoted by 
\begin{equation}
w_{i} =   (a_{1i},  \cdots , a_{ni})^{T}         
\in \mathbb{Q}^{n}
\nonumber
\end{equation}

\noindent
Consider the relation $c_{1}^{2}-2c_{2}+p_{1}=0$ ,  and apply   $p^{*} \circ h$ to it. 
\begin{equation}
p^{*}\circ h(c_{1}^{2}-2c_{2}+p_{1})= \widetilde{h}\circ p^{*}(c_{1}^{2}-2c_{2}+p_{1}) = \widetilde{h}(t_{1}^{2}+...+t_{n}^{2})=0  
\nonumber
\end{equation}
Now,  
\begin{equation}
 \widetilde{h}(t_{1}^{2}+...+t_{n}^{2})= \Sigma_{i=1}^{n}\langle w_{i}, w_{i}\rangle t_{i}^{2}+ 2 \Sigma_{i<j} \langle w_{i}, w_{j}\rangle t_{i}t_{j}=0 
 \label{eqn621}
\end{equation}
Since , $t_{1}^{2}+...+t_{n}^{2}=0$ ,  from \ref{eqn621},  we get
$ D :=  \langle w_{i}, w_{i} \rangle =...= \langle w_{j},w_{j} \rangle  \geq 0$
for all $i,  j \in \{1,  \cdots ,  n \}$,  and  
$\langle w_{i}, w_{j} \rangle = 0$
for all $i\neq j$ ,  $i,  j \in \{ 1,2,...,n\}$. \\
 \noindent 
 Let,  if possible,  $D \neq 0$.  Then $H$ is invertible,  which is a contradiction.  Therefore,  $D=0$,  and $H$ is the zero matrix,  and $h$ is the zero endomorphism.

\end{proof}

\noindent
We have an analogous result,  when the image of the first Pontrjagin class of $\xi$ vanishes.  
 \begin{myprop}   {\bf [Theorem \ref{prop64}]}  
 Let $h:H^{*}(M_{n,k}; \Q) \to H^{*}(M_{n,k}; \Q)$  be an endomorphism of the cohomology algebra,  such that $h(p_{1})=0$,  where $p_{1}$ denotes the first Pontrjagin class of $\xi$.  If $h$ admits a splitting $\widetilde{h}$,   then $\widetilde{h} =0$.  In such a case,  $h$ will be the zero endomorphism of $H^{*}(M_{n,k};\Q)$.\\
\end{myprop}

\begin{proof}
Recall the following commuting diagram: \\
 \begin{tikzcd}
H^{*}(M_{n,k} ; \Q)\arrow[r,"h"] \arrow[d,"p^{*}"] & H^{*}(M_{n,k} ; \Q) \arrow[d ,"p^{*}"] \\ 
 H^{*}(SO(2n)/T^{n}; \Q ) \arrow[r,"\widetilde{h}"] & H^{*}(SO(2n)/T^{n}; \Q)
\end{tikzcd}
\\
Consider the matrix representative of the endomorphism $\widetilde{h}$,   $H = (a_{ij})_{1 \leq i, j \leq n}$.  \\
\noindent
Now,  $h(p_{1})=0$ implies $p^{*}\circ h(p_{1})= \widetilde{h}\circ p^{*}(p_{1})=\widetilde{h}(t_{k+1}^{2}+...+t_{n}^{2})=0$.\\
\noindent
Let 
\begin{equation}
u_{i} =    (a_{1i} ,  \cdots ,  a_{ki})^{T}  
\in \mathbb{Q}^{k} ; \hspace{0.5cm}  v_{i} = (a_{k+1,i},  \cdots ,  a_{n,i})^{T} \in \mathbb{Q}^{n-k}           
\nonumber
\end{equation}
for $i \in \{ 1,  \cdots,n \}$.  \\
Since 
$\widetilde{h}(t_{k+1}^{2}+ \cdots +t_{n}^{2})= \Sigma_{i=1}^{n}\langle v_{i},v_{i}\rangle t_{i}^{2}+ 2\Sigma_{ i<j} \langle v_{i} , v_{j} \rangle t_{i}t_{j} =0$
and $t_{1}^{2}+ \cdots +t_{n}^{2}=0$ ,  we get $D =  \langle v_{1}, v_{1} \rangle = \cdots = \langle v_{n},v_{n}\rangle \geq 0 $ and $ \langle v_{i}, v_{j} \rangle = 0$,  for $ i, j \in \{ 1,...,n\},   i \neq j$. \\
   
\noindent
Let,  if possible,  $D \neq 0$.  Then we have $n$ linearly independent vectors in $\Q^{n-k}$ ,  which is not possible.  Therefore,  $D=0$,  and $v_{1}= \cdots =v_{n}=\overline{0} \in \Q^{n-k}$. \\

\noindent
So the matrix $H$ is of the following form,  with the last $n-k$ rows being zero:  

\begin{equation}
H=
\begin{bmatrix}
 a_{11} &   \cdots & a_{1k} & a_{1 k+1} & \cdots & a_{1n} \\
 \cdots &   \cdots& \cdots & \cdots & \cdots& \cdots \\
 a_{k1}  &    \cdots &  a_{kk} &  a_{k k+1} & \cdots & a_{kn} \\
  0  &    \cdots &  0 &  0 & \cdots & 0 \\
 \cdots &   \cdots& \cdots & \cdots & \cdots& \cdots  \\
  0 &   \cdots & 0 & 0 & \cdots & 0
\end{bmatrix}
\nonumber 
\end{equation}

\noindent
Using the relation $c_{1}^{2}-2c_{2}+p_{1}=0$,  we get, 

$p^{*}\circ h(c_{1}^{2}-2c_{2}+p_{1})= \widetilde{h}\circ p^{*}(c_{1}^{2}-2c_{2}+p_{1})= \widetilde{h}(t_{1}^{2}+ \cdots +t_{n}^{2})=0$ \\
\noindent
We have  $\widetilde{h}(t_{1}^{2}+ \cdots +t_{n}^{2})= \Sigma_{i=1}^{n}\langle u_{i},u_{i}\rangle t_{i}^{2} + 2 \Sigma_{i\neq j , i< j}\langle u_{i}, u_{j}\rangle t_{i}t_{j}=0$
and $t_{1}^{2}+ \cdots +t_{n}^{2}=0$.  
Therefore,  $D_{1} :=   \langle u_{1}, u_{1} \rangle = \cdots = \langle u_{n},u_{n}\rangle  \geq 0$ 
and $\langle u_{i}, u_{j} \rangle = 0$,   where $i \neq j \in \{ 1,...,n\}.$

\noindent
Let,  if possible,  $D_{1} \neq 0$.  Then,  we get $n$ linearly independent vectors in $\Q^{k}$,  which is not possible.  Hence $D_{1}=0$ and $u_{1}= \cdots =u_{n}=\overline{0} \in \Q^{k}$.  Therefore,  the matrix $H$ is the zero matrix,  and the endomorphism $h$ is the zero endomorphism. \\
\end{proof}

\begin{myprop}  {\bf [Theorem \ref{prop65}]}  
Let $h:H^{*}(M_{n,k}; \Q) \to H^{*}(M_{n,k}; \Q)$  be an endomorphism  such that $h(e_{2n-2k})=0$ , where $n-k=2$ and $e_{2n-2k}$ denotes the Euler class of $\xi$.   If $h$ admits a splitting $\widetilde{h}$,   then $\widetilde{h} =0$.  and $h$ is also the zero endomorphism. 
\end{myprop}

\begin{proof}
Since $n-k =2$,  we have $h(e_{4})=0$,  and 
$p^{*}\circ h(e_{4})= \widetilde{h}\circ p^{*}(e_{4})= \widetilde{h}(t_{k+1}t_{k+2})=0$.\\

\noindent
The matrix $H$ is  of the form : 
\begin{equation}
\begin{bmatrix}
 a_{11} & a_{12} & \cdots & a_{1k} & a_{1k+1} &  a_{1k+2} \\
 \cdots &  \cdots & \cdots & \cdots &  \cdots & \cdots &   \\
 a_{k1}  &  a_{k2} & \cdots &  a_{kk} &  a_{kk} & a_{kk+2} \\
  b_{1}  &  b_{2} & \cdots &  b_{k} &   d_{1}  &   d_{2} \\
  f_{1}  &  f_{2} & \cdots &  f_{k} &   d_{3} &   d_{4}
\end{bmatrix} 
\end{equation}

\noindent
Now,  $0= \widetilde{h}(t_{k+1}t_{k+2}) 
= \Sigma_{i=1}^{k} (b_{i}f_{i})t_{i}^{2}+ (d_{1}d_{3})t_{k+1}^{2}+(d_{2}d_{4})t_{k+2}^{2}+\Sigma_{i\neq j \\ i,j \in \{1,...,k\}}(b_{i}f_{j}+b_{j}f_{i}) t_{i}t_{j}+ (d_{1}d_{4}+d_{2}d_{3})t_{k+1}t_{k+2}+ \Sigma_{i=1}^{k}(b_{i}d_{3}+f_{i}d_{1})t_{i}t_{k+1}+\Sigma_{i=1}^{k}(b_{i}d_{4}+f_{i}d_{2})t_{i}t_{k+2}$. \\

\noindent
Comparing this with $t_{1}^{2}+...+t_{k}^{2}+t_{k+1}^{2}+t_{k+2}^{2}=0$ ,  we get the following
\begin{equation}
    b_{1}f_{1}= \cdots =b_{k}f_{k}=d_{1}d_{3}=d_{2}d_{4}
    \label{eqn624}
\end{equation}
\begin{equation}
    b_{i}f_{j}+b_{j}f_{j}=0 ,  i\neq j \in \{1,  \cdots ,k\}
    \nonumber
\end{equation}
\begin{equation}
    d_{1}d_{4}+d_{2}d_{3}=0 
    \label{eqn625}
\end{equation}
\begin{equation}
    (b_{i}d_{3}+f_{i}d_{1})=0 ; \ (b_{i}d_{4}+f_{i}d_{2})=0
    \nonumber
\end{equation}
 for all $ i \in \{1, 2,  \cdots,  k \}$.  \\
 
\noindent
From \ref{eqn624} and \ref{eqn625},  we get  $d_{1}d_{3}=d_{2}d_{4}$.  Let $d_{1}d_{3}=d_{2}d_{4}\neq 0$ and hence,   $d_{i}\neq 0$ for all $i$.  Therefore we have,  $\dfrac{d_{2}d_{4}}{d_{3}}d_{4}+d_{2}d_{3}=0$,  that is,   
    $\dfrac{d_{2}}{d_{3}}(d_{4}^{2}+d_{3}^{2})=0$. 
This implies $d_{3}=d_{4}=0$,  which is a contradiction.  Therefore,  we have $d_{1}d_{3}=d_{2}d_{4}=0$.\\

\noindent
From \ref{eqn624} , we get 
\begin{equation}
    b_{1}f_{1}= \cdots = b_{k}f_{k}=d_{1}d_{3}=d_{2}d_{4}=0
    \nonumber
\end{equation}

\noindent
 A similar computation for $h(p_{1})$ gives us: 
\begin{equation}
    b_{1}^{2}+ f_{1}^{2}=b_{2}^{2}+f_{2}^{2}= \cdots =b_{k}^{2}+f_{k}^{2}
    \nonumber
\end{equation}
\begin{equation}
    d_{1}^{2}+d_{3}^{2}=d_{2}^{2}+d_{4}^{2}
    \nonumber
\end{equation}
\begin{equation}
    d_{1}b_{i}+d_{3}f_{i}= d_{2}b_{i}+d_{4}f_{i}=0 , \ i \in \{ 1,2,  \cdots ,k\}
    \nonumber
\end{equation}

\noindent
Analyzing the above equations we observe that the last two rows of the matrix $H$ are of one of the following forms: \\

$1) \begin{bmatrix}
 
  0  &  0 & \cdots &  0 &  0 & 0 \\
  0  &  0 & \cdots &  0 &  0 & 0
  
\end{bmatrix} $ \\

$2) \begin{bmatrix}
 
  0  &  0 & \cdots &  0 &  0 & 0 \\
  f_{1}  &  f_{1} & \cdots &  f_{1} &  0 & 0
  
\end{bmatrix} $ \\

$3) \begin{bmatrix}
 
  b_{1}  &  b_{1} &  \cdots &  b_{1} &  0 & 0 \\
  0  &  0 & \cdots &  0 &  0 & 0
  
\end{bmatrix}$\\

$4) \begin{bmatrix}
 
  0  &  0 & \cdots &  0 &  0 & 0 \\
  0  &  0 &  \cdots &  0 &  d_{3} & \pm d_{3}
  
\end{bmatrix}$ \\

$5) \begin{bmatrix}
 
  0  &  0 & \cdots &  0 &  d_{1} & \pm d_{1} \\
  0  &  0 & \cdots &  0 &  0 & 0
  
\end{bmatrix}$ \\

\noindent
Therefore,  the matrix $H$ has at least one row with all zeros,  and hence $H$ is singular.  We have already proved that if $H$ is singular,  then $H$ is the zero matrix.  (Note that this statement is valid for $\Q$ and $\R$.) Hence,  $h$  is the zero endomorphism.

\end{proof}

\vspace{0.5cm}

\section{Lefschetz Number of Certain Maps}

\noindent
Let $N = \dim_{\C} M_{n,k}$ and let $d_{2j} = \dim H^{2j} (M_{n,k} ; \Q)$.  
Then by Poincare Duality,  $d_{2j} = d_{2N - 2j}$.  
The Lefschetz number,  $L(f)$, of a map $f \colon M_{n,k} \longrightarrow M_{n,k}$ is an invariant, connected to the fixed points of the map.  It is given by : 
\begin{displaymath}
L(f) = \sum_{i=0}^{2N} (-1)^{i} \Tr (f_{i}^{*} \colon H^{i} (M_{n,k} ; \Q) \longrightarrow H^{i} (M_{n,k} ; \Q)) =  \sum_{i=0}^{2N} (-1)^{i} \Tr (f_{i}^{*})
\end{displaymath}
where $\Tr$ denotes trace.  A version of Lefschetz Fixed Point Theorem states that if $L(f) \neq 0$,  then $f$ has at least one fixed point.  \\

\noindent
Let $m \neq 0$ and $f \colon H^{*} (M_{n,k} ; \Q) \longrightarrow H^{*} (M_{n,k} ; \Q)$ be an endomorphism of the cohomology algebra such that : 
\begin{eqnarray}
    f(c_{i}) &=& m^{i} c_{i} ; \hspace{0.5cm} 1 \leq i \leq k \nonumber \\ 
    f(p_{j}) &=& m^{2j} p_{j} ;  \hspace{0.5cm} 1 \leq j \leq n-k  \nonumber \\
    f(e_{2n-2k}) &=& \pm m^{n-k}e_{2n-2k} \nonumber
\end{eqnarray}

\noindent
 The Lefschetz number, $L(f)$,  for the map $f$,  is defined to be : 
 \begin{displaymath}
  L(f) :=  \sum_{i=0}^{2N} (-1)^{i} \Tr (f_{i})  
  = 1\cdot d_{0} + md_{2} + m^{2} d_{4} + \cdots + m^{N} d_{2N} 
  \end{displaymath}
  \begin{displaymath}
  = (1+m^{N}) d_{0} + (m + m^{N-1}) d_{2} + (m^{2} + m^{N-2}) d_{4} + \cdots 
\end{displaymath}

\noindent
Then we have the following results: 
\begin{prop}
Let $m \neq 0$  be an integer,  and $f \colon H^{*} (M_{n,k} ; \Q) \longrightarrow H^{*} (M_{n,k} ; \Q)$ be an endomorphism of the cohomology algebra such that : 
\begin{eqnarray}
    f(c_{i}) &=& m^{i} c_{i} ; \hspace{0.5cm} 1 \leq i \leq k \nonumber \\ 
    f(p_{j}) &=& m^{2j} p_{j} ;  \hspace{0.5cm} 1 \leq j \leq n-k  \nonumber \\
    f(e_{2n-2k}) &=& m^{n-k}e_{2n-2k} \nonumber
\end{eqnarray}

\noindent
\begin{enumerate}
\item [(A)] If $m > 0$,  then $L(f)$ is positive.  If,  in addition,  $m$ is coprime to the order of the Weyl group of $SO(2n)$,  then $f$ is induced by a self-map $F \colon M_{n,k} \longrightarrow M_{n,k}$ ($F^{*} \simeq f$).  In particular,  this self-map $F$ will have a fixed point.  \\

\noindent
\item[(B)] If $m = -1$,  then,  
\begin{itemize}
\item[(i)] if $N$ is odd,  then $L(f) =0$ 
\item[(ii)] if $N$ is even,  then 
\begin{displaymath}
L(f) = 2 \sum_{j=0}^{\frac{N-2}{2}} (-1)^{j} d_{2j}  + (-1)^{\frac{N}{2}} d_{N}
\end{displaymath}
\end{itemize}
\end{enumerate}
\end{prop}

\begin{proof}

\noindent
 \begin{displaymath}
  L(f) :=  \sum_{i=0}^{2N} (-1)^{i} \Tr (f_{i})  
  = (1+m^{N}) d_{0} + (m + m^{N-1}) d_{2} + (m^{2} + m^{N-2}) d_{4} + \cdots 
\end{displaymath}
Therefore,  $m>0$ implies $L(f) >0$.  
Note that when $m$ is an integer and $ f(e_{2n-2k}) =  m^{n-k}e_{2n-2k}$,  $f$ is an Adams map.  Further,  if $m$ is coprime to the order of the Weyl group of $SO(2n)$,  then by Lin's result  (Theorem 1.3,  \cite{lin}),  $f$ is induced by a self-map $F \colon M_{n,k} \longrightarrow M_{n,k}$ (that is,  $F^{*} \simeq f)$.  Since the Lefschetz number of $F$ is positive,  by the Lefschetz Fixed Point Theorem,  $F$ has a fixed point.  Hence we have part (A).\\

\noindent
We have the following two cases: \\
\noindent
{\bf Case:} When $N$ is odd,  
\begin{displaymath}
L(f) = (1+m^{N}) d_{0} + (m + m^{N-1}) d_{2} + (m^{2} + m^{N-2}) d_{4} + \cdots +  (m^{\frac{N-1}{2}} + m^{\frac{N+1}{2}}) d_{N-1}
\end{displaymath}
\begin{displaymath}
= \sum_{j=0}^{\frac{N-1}{2}} (m^{j} + m^{N-j}) d_{2j}
\end{displaymath}

\noindent
{\bf Case:} When $N$ is even,  
\begin{displaymath}
L(f) =  (1+m^{N}) d_{0} + (m + m^{N-1}) d_{2} + (m^{2} + m^{N-2}) d_{4} + \cdots +  (m^{\frac{N-2}{2}} + m^{\frac{N+2}{2}}) d_{N-2} + m^{\frac{N}{2}} d_{N}
\end{displaymath}
\begin{displaymath}
= \sum_{j=0}^{\frac{N-2}{2}} (m^{j} + m^{N-j}) d_{2j} + m^{\frac{N}{2}} d_{N}
\end{displaymath}
\noindent
Putting $m=-1$,  we get part (B) of the proposition.  
\end{proof}

\vspace{0.3cm}

\begin{prop}
Let $m$  be an integer,  and $f \colon H^{*} (M_{n,k} ; \Q) \longrightarrow H^{*} (M_{n,k} ; \Q)$ be an endomorphism of the cohomology algebra such that : 
\begin{eqnarray}
    f(c_{i}) &=& m^{i} c_{i} ; \hspace{0.5cm} 1 \leq i \leq k \nonumber \\ 
    f(p_{j}) &=& m^{2j} p_{j} ;  \hspace{0.5cm} 1 \leq j \leq n-k  \nonumber \\
    f(e_{2n-2k}) &=& m^{n-k}e_{2n-2k} \nonumber
\end{eqnarray}
If $m <0$,  $m \neq -1$,  and $N$ and $\dfrac{N-1}{2}$ are odd,  and $|m|$ is greater than all the Betti numbers,  then $L(f)$ is negative.  
\end{prop}

\begin{proof}
\noindent
  When $N$ is odd ,  the Lefschetz number is given by 
   \begin{equation}
       L(f)= \sum_{j=0}^{\frac{N-1}{2}}(m^{j}+m^{N-j})d_{2j} \nonumber
   \end{equation}
 
 \noindent
Substituting  $-m$ in place of $m$,  we have 
\begin{eqnarray}
       L(f)&=& \sum_{j=0}^{\frac{N-1}{2}}\left [(-m)^{j}+(-m)^{N-j} \right ]d_{2j}  \nonumber \\
           &=& - \sum_{j=0}^{\frac{N-3}{4}} ( m^{N-2j}-m^{2j})d_{4j}+ \sum_{j=0}^{\frac{N-3}{4}} ( m^{N-(2j+1)}-m^{2j+1})d_{4j+2}  \nonumber \\ 
           &=& -B + A  \hspace{0.3cm} (say) \nonumber
   \end{eqnarray}

\noindent
Let $ M_{1} = \max \{ d_{2}, d_{6}, ... , d_{N-1}\} $ ,   $ m_{1} =\min \{ d_{2}, d_{6},...,d_{N-1}\} $ ,  $ M_{2} = \max \{ d_{0}, d_{4}, ... , d_{N-3}\} $,  and $ m_{2} = \min \{ d_{0}, d_{4},  \cdots,  d_{N-3}\} $. 
Then  
\begin{equation}
  m_{1} \left [\sum_{j=0}^{\frac{N-3}{4}} ( m^{N-(2j+1)}-m^{2j+1}) \right ] \leq A \leq M_{1} \left [\sum_{j=0}^{\frac{N-3}{4}} ( m^{N-(2j+1)}-m^{2j+1}) \right ] \nonumber
\end{equation}
and 
\begin{equation}
  m_{2} \left [ \sum_{j=0}^{\frac{N-3}{4}} ( m^{N-2j}-m^{2j}) \right ] \leq B \leq M_{2} \left [ \sum_{j=0}^{\frac{N-3}{4}} ( m^{N-2j}-m^{2j})\right ] \nonumber
\end{equation}
\noindent
Combining the above inequalities,  we get, 
\begin{equation}
    m_{1} \left ( \sum_{j=0}^{\frac{N-3}{4}} ( m^{N-(2j+1)}-m^{2j+1}) \right )-  M_{2} \left ( \sum_{j=0}^{\frac{N-3}{4}} ( m^{N-2j}-m^{2j}) \right ) \leq (A-B) = L(f) \label{eqn728}
\end{equation}
and 
\begin{equation}
   L(f)=(A-B)
   \leq M_{1} \left ( \sum_{j=0}^{\frac{N-3}{4}} ( m^{N-(2j+1)}-m^{2j+1}) \right ) - m_{2} \left ( \sum_{j=0}^{\frac{N-3}{4}} ( m^{N-2j}-m^{2j}) \right )  \label{eqn729}
\end{equation}
 \noindent
 Let $m \neq 1$ and $p, q$  be positive integers and $k,l$ be non-negative integers,  such that $k-l=2$,  and $p+1 = q+k$,   and $p = q+l+1$.  Then we have the inequality $ m^{p+1}-m^{q} > m^{p}-m^{q+1}$.  Putting $p = N-1$ and $q=0$,  we get $(m^{N}-1) > (m^{N-1}-m)$.  This implies 
\begin{equation}
  \sum_{j=0}^{\frac{N-3}{4}} ( m^{N-2j}-m^{2j}) > \sum_{j=0}^{\frac{N-3}{4}} ( m^{N-(2j+1)}-m^{2j+1})  
\end{equation}
\noindent
For any $m_{1}$ and $M_{2}$,  the LHS of Equation (\ref{eqn728}) is always negative.  Since $m_{2} =1$,  the RHS of Equation (\ref{eqn729}) is negative if 
\begin{equation}
    M_{1}<< \frac{ \left (\sum_{j=0}^{\frac{N-3}{4}} ( m^{N-2j}-m^{2j}) \right )}{\left ( \sum_{j=0}^{\frac{N-3}{4}} ( m^{N-(2j+1)}-m^{2j+1}) \right )} \sim m \nonumber
\end{equation}
\noindent
Therefore,  $M_{1} << m $ implies $L(f)$ is negative.  
\end{proof}

\vspace{0.5cm}

\vspace{0.5cm}

\end{document}